\documentclass[12pt,reqno,final]{amsart}

\usepackage{amsmath}
\usepackage{amsthm}
\usepackage{amsfonts}  
\usepackage{amssymb}   
\usepackage{latexsym}

\theoremstyle{plain}
\newtheorem{theorem}{Theorem}
\newtheorem{proposition}[theorem]{Proposition}
\newtheorem{corollary}[theorem]{Corollary}

\theoremstyle{remark}
\newtheorem{remark}{Remark}

\DeclareMathAlphabet\mathoo{U}{eur}{b}{n}
\DeclareMathOperator{\Real}{Re}

\usepackage[dvips]{graphicx}
\usepackage[hypertex,colorlinks=true,citecolor=blue]{hyperref} 

\begin{document}
\title[Computation of Green's function]{Computation of Green's function\\ of the bounded solutions problem}

\author{V.G. Kurbatov}
 \email{kv51@inbox.ru}
 \address{Department of Mathematical Physics,
Voronezh State University\\ 1, Universitetskaya Square, Voronezh 394018, Russia}

\author{I.V. Kurbatova}
 \email{la\_soleil@bk.ru}
 \address{Department of Software Development and Information Systems Administration,
Voronezh State University\\ 1, Universitetskaya Square, Voronezh 394018, Russia}



\begin{abstract}
It is well known that the equation $x'(t)=Ax(t)+f(t)$, where $A$ is a square matrix, has a unique bounded solution $x$ for any bounded continuous free term $f$, provided the coefficient $A$ has no eigenvalues on the imaginary axis. This solution can be represented in the form
\begin{equation*}
x(t)=\int_{-\infty}^{\infty}\mathcal G(t-s)x(s)\,ds.
\end{equation*}
The kernel $\mathcal G$ is called Green's function. In the paper, a representation of Green's function in the form of the Newton interpolating polynomial is used for approximate calculation of $\mathcal G$. An estimate of the sensitivity of the problem is given.
\end{abstract}

\maketitle

\section*{Introduction}\label{s:Introduction}
An approximate computation of analytic functions $f$ of matrices is an important problem in modern numerical mathematics~\cite{Frommer-Simoncini,Golub-Van_Loan96,Higham08,Higham-Al-Mohy2010,Ikramov1991a:eng}. The class of specific functions $f$ that are usually discussed is not very large. These are the exponential function, the cosine and the sine (and some of their modifications), the square root and other non-integer powers, the logarithm, and the sign function.
In this paper, we want to call attention to the function $g_t$ arising in calculating Green's function (see formula~\eqref{e:Green's function}) for the problem of bounded on the axis solutions of the differential equation $y'(t)-Ay(t)=f(t)$.

If the eigenvalues of $A$ are arranged in a special order, then the Newton interpolating polynomial of $g_t$ takes a simple form (Theorem~\ref{t:p_t}). It allows us to reduce the problem of calculating $g_t(A)$ to the substitution of $A$ into the interpolating polynomial of the functions $\widetilde{\exp^\pm_t}$ which degree is a half of the degree of $g_t$.

The norm of the differential of Green's function with respect to $A$ is called the condition number. It shows how inaccuracies in $A$ influence on $g_t(A)$. We derive an estimate of the condition number from above (Corollary~\ref{c:est of dg_t}).

In Sections~\ref{s:Matrix functions}, \ref{s:divided differences}, and~\ref{s:Newton polynomial}, preliminaries are collected. In Section~\ref{s:Green}, we recall the definition of Green's function and some its properties and describe its representation in the form of the Newton interpolating polynomial. In Section~\ref{s:sensitivity} we discuss the sensitivity of the problem under consideration. In Section~\ref{s:algorithm} we present our algorithm. The algorithm is convenient for symbolic calculations. In Section~\ref{s:experiments} we describe some numerical experiments.

\section{Matrix functions}\label{s:Matrix functions}

Let $A$ be a complex $N\times N$-matrix. Let $\mathbf1$ be the identity matrix. The polynomial
\begin{equation*}
p_A(\lambda)=\det(\lambda\mathbf1-A)
\end{equation*}
is called the \emph{characteristic polynomial} of the matrix $A$. Let $\lambda_1$, \dots, $\lambda_M$ be the complete set of the roots of the characteristic polynomial $p_A$, and $n_1$, \dots, $n_M$ be their multiplicities; thus $n_1+\dots+n_M=N$. It is well known that $\lambda_1$, \dots, $\lambda_M$ are \emph{eigenvalues} of $A$. The numbers $n_k$ are called (\emph{algebraic{\rm)} multiplicities} of the eigenvalues $\lambda_k$. The set $\sigma(A)=\{\,\lambda_1, \dots, \lambda_M\,\}$ is called the \emph{spectrum} of $A$.

Let $U\subseteq\mathbb C$ be an open set that contains the spectrum $\sigma(A)$.
Let $f:\,U\to\mathbb C$ be an analytic function. The \emph{function $f$ of the matrix} $A$ is defined~\cite[ch.~V, \S~1]{Hille-Phillips:eng},~\cite[p.~17]{Daletskii-Krein:eng} by the formula
\begin{equation*}
f(A)=\frac1{2\pi i}\int_\Gamma f(\lambda)(\lambda\mathbf1-A)^{-1})\,d\lambda,
\end{equation*}
where the contour $\Gamma$ surrounds the spectrum $\sigma(A)$.

\begin{proposition}[{\rm\cite[Theorem 5.2.5]{Hille-Phillips:eng}}]\label{p:func calc}
The mapping $f\mapsto f(A)$ preserves algebraic operations, i.~e.,
\begin{align*}
(f+g)(A)&=f(A)+g(A),\\
(\alpha f)(A)&=\alpha f(A),\\
(fg)(A)&=f(A)g(A),
\end{align*}
where $f+g$, $\alpha f$ and $fg$ are defined pointwise.
\end{proposition}

\begin{proposition}[{\rm see, e.g., \cite[Proposition 2.3]{Frommer-Simoncini}}]\label{p:f=p}
Let functions $f$ and $g$ be analytic in a neighbourhood of the spectrum $\sigma(A)$ and
\begin{equation*}
f^{(j)}(\lambda_k)=g^{(j)}(\lambda_k),\qquad k=1,\dots,M;\;j=0,1,\dots,n_k-1.
\end{equation*}
Then
\begin{equation*}
f(A)=g(A).
\end{equation*}
\end{proposition}

\section{Divided differences}\label{s:divided differences}

Let $\mu_1$, $\mu_2$, \dots, $\mu_N$ be given complex numbers (some of them may coincide
with others) called \emph{points of interpolation}. Let a complex-valued function $f$ be
defined and analytic in a neighbourhood of these points. \emph{Divided differences}
of the function $f$ with respect to the points $\mu_1$,
$\mu_2$, \dots, $\mu_N$ are defined (see, e.g.,~\cite{
Gelfond:eng,Jordan}) by the recurrent
relations
 \begin{equation}\label{e:divided differences}
 \begin{split}
f[\mu_i]&=f(\mu_i),\\
f[\mu_i,\mu_{i+1}]&=\frac{f(\mu_{i+1})-f(\mu_i)}{\mu_{i+1}-\mu_i},\\
f[\mu_i,\dots,\mu_{i+m}]&=\frac{f(\mu_{i+1},\dots,\mu_{i+m})
-f(\mu_{i},\dots,\mu_{i+m-1})} {\mu_{i+m}-\mu_i}.
 \end{split}
 \end{equation}
In these formulas, if the denominator vanishes, then the quotient is understood as the derivative with respect to one of the arguments of the previous divided difference (this agreement can by derived by continuity from Corollary~\ref{c:f[] is continuous}).

\begin{proposition}[{\rm\cite[ch.~1, formula (54)]{Gelfond:eng}}]\label{p:f[] via Gamma}
Let the function $f$ be analytic in a neighbourhood of the points of interpolation $\mu_1$, $\mu_2$, \dots, $\mu_N$. Then divided differences possess the representation
\begin{equation*}
f[\mu_{1},\dots,\mu_{N}]=\frac1{2\pi i}\int_{\Gamma}\frac{f(z)}{\Omega(z)}\,dz,
\end{equation*}
where the contour $\Gamma$ encloses all the points of interpolation and
\begin{equation*}
\Omega(z)=\prod_{k=1}^N(z-\mu_k).
\end{equation*}
\end{proposition}

\begin{corollary}\label{c:f[] is continuous}
Divided differences are differentiable functions of their arguments.
\end{corollary}
\begin{proof}
The statement follows from Proposition~\ref{p:f[] via Gamma}.
\end{proof}

 \begin{corollary}\label{c:f[] is indep of order}
Divided differences $f[\mu_{1},\dots,\mu_{N}]$ are
symmetric function, i.e., they do not depend on the order of their arguments $\mu_1$, \dots, $\mu_{N}$.
 \end{corollary}
 \begin{proof}
The statement follows from Proposition~\ref{p:f[] via Gamma}.
 \end{proof}

\begin{proposition}[{\rm\cite[ch.~1, formula (48)]{Gelfond:eng}}]\label{p:repr of Delta}
Let the points of interpolation $\mu_j$ be distinct. Then divided differences possess the representation
\begin{equation*}
f[\mu_{1},\dots,\mu_{N}]=\sum_{j=1}^N\frac{f(\mu_j)}{\prod\limits_{\substack{k=1\\k\neq j}}^N(\mu_j-\mu_k)}.
\end{equation*}
\end{proposition}
 \begin{proof}
The statement follows from Proposition~\ref{p:f[] via Gamma}.
 \end{proof}

\section{The Newton interpolating polynomial}\label{s:Newton polynomial}

The set $\lambda_1,\dots,\lambda_M\in\mathbb C$ of \emph{interpolation points} together with the set $n_1,\dots,n_M\in\mathbb N$ of their \emph{multiplicities} is called \emph{multiple interpolation data}.
We set $N=n_1+\dots+n_M$.

Let $U\subseteq\mathbb C$ be an open neighbourhood of the set $\lambda_1,\dots,\lambda_M$ of the points of interpolation and $f:\,U\to\mathbb C$ be an analytic function.
An \emph{interpolating polynomial} of $f$ that corresponds to the multiple interpolation data is a polynomial $p$ of degree степени $N-1$ satisfying the equalities
\begin{equation*}
p^{(j)}(\lambda_k)=g^{(j)}(\lambda_k),\qquad k=1,\dots,M;\;j=0,1,\dots,n_k-1.
\end{equation*}

\begin{proposition}\label{p:f->p}
Let $A$ be a complex $N\times N$-matrix.
Let $p$ be an interpolating polynomial of $f$ that corresponds to the points $\lambda_1,\dots,\lambda_M$ of the spectra of the matrix $A$ counted according to their
multiplicities $n_1,\dots,n_M$.
Then
\begin{equation*}
p(A)=f(A).
\end{equation*}
\end{proposition}
\begin{proof}
The statement immediately follows from Proposition~\ref{p:f=p}.
\end{proof}

\begin{remark}\label{r:n_k}
Propositions~\ref{p:f=p} and~\ref{p:f->p} remain valid if one assumes that $n_1$, \dots, $n_M$ are the maximal sizes of the corresponding Jordan blocks. Formally, this assumption decreases the degree $N-1$ of the interpolating polynomial. But in approximate calculations, it does not help. In fact, the property of coincidence of eigenvalues and the sizes of Jordan blocks are numerically unstable. Thus in practice, one can hardly meet Jordan blocks of the size more than $1\times1$.
\end{remark}

\begin{proposition}[{\rm\cite[p.~20]{Jordan}}]\label{p:Newton poly}
For any analytic function $f$, the interpolation polynomial exists and unique.
Let $\mu_1,\dots,\mu_N$ be the points of multiple interpolation data $\lambda_1,\dots,\lambda_M$, listed in an arbitrary order and repeated as many times as their multiplicities $n_1,\dots,n_M$.
Then the interpolating polynomial possesses the representation
\begin{equation}\label{e:poly Newton}
\begin{split}
p(z)&=f[\mu_1]+f[\mu_1,\mu_2](z-\mu_1)+
f[\mu_1,\mu_2,\mu_3](z-\mu_1)(z-\mu_2)\\
&+f[\mu_1,\mu_2,\mu_3,\mu_4](z-\mu_1)(z-\mu_2)(z-\mu_3)+\dots\\
&+f[\mu_1,\mu_2,\dots,\mu_N](z-\mu_1)(z-\mu_2)\dots(z-\mu_{N-1}).
 \end{split}
\end{equation}
\end{proposition}

Representation~\eqref{e:poly Newton} is called~\cite{Gelfond:eng,Jordan} the \emph{interpolating polynomial in the Newton form} or shortly \emph{the Newton
interpolating polynomial} with respect to the points $\mu_1$, $\mu_2$, \dots, $\mu_N$.

\section{Green's function}\label{s:Green}
In this Section, we recall the definition and prove some properties of Green's function.

Let $A$ be a complex $N\times N$-matrix. We consider the differential equation
\begin{equation}\label{e:x'=Ax+f}
x'(t)=Ax(t)+f(t),\qquad t\in\mathbb R.
\end{equation}
We are interested in \emph{bounded solutions problem}, i.e. seeking bounded solution $x:\,\mathbb R\to\mathbb C^N$ under the assumption that the free term $f:\,\mathbb R\to\mathbb C^N$ is a bounded function. The bounded solutions problem has its origin in the work of Perron~\cite{Perron}. Its different modifications can be found in~\cite{BaskakovMS15:eng,Daletskii-Krein:eng,Henry81:eng,KurbatovSMZ86:eng,Pechkurov15:eng}. See also references therein.

Suppose that $\sigma(A)$ does not intersect the imaginary axis. In this case the functions
\begin{align*}
\exp^+_t(\lambda)&=\begin{cases}
e^{\lambda t}, & \text{if $\Real\lambda<0$},\\
0, & \text{if $\Real\lambda>0$},\end{cases}\\
\exp^-_t(\lambda)&=\begin{cases}
0, & \text{if $\Real\lambda<0$},\\
e^{\lambda t}, & \text{if $\Real\lambda>0$},\end{cases}\\
g_t(\lambda)&=\begin{cases}
-\exp^-_t(\lambda), & \text{if $t<0$},\\
\exp^+_t(\lambda), & \text{if $t>0$}
\end{cases}
\end{align*}
are analytic in a neighbourhood of the spectrum $\sigma(A)$. We set
\begin{equation}\label{e:Green's function}
\mathcal G(t)=g_t(A),\qquad t\neq0.
\end{equation}
The function $\mathcal G$ is called~\cite{Daletskii-Krein:eng} \emph{Green's function} of the boundary solutions problem for equation~\eqref{e:x'=Ax+f}.

The following proposition is well known.

\begin{proposition}\label{p:G-properties}
Green's function possesses the properties{\rm:}
\begin{enumerate}
 \item $P^+=\mathcal G(+0)$ and $P^-=\mathcal G(-0)$ are projectors, i.~e. they satisfy the identity $P^2=P$,
 \item $P^+-P^-=\mathbf1$,
 \item $\mathcal G(t_1)\mathcal G(t_2)=\mathcal G(t_1+t_2)$,
 \item $\mathcal G(t_1)\mathcal G(t_2)=0$ for $t_1t_2<0$.
 \item $\frac d{dt}\mathcal G(t)=A\mathcal G(t)$ for $t\neq0$.
\end{enumerate}
\end{proposition}
\begin{proof}
The statement follows from Proposition~\ref{p:func calc} and the identities
\begin{gather*}
g_{\pm0}^2(\lambda)=g_{\pm0}(\lambda),\\
g_{+0}(\lambda)-g_{-0}(\lambda)=1,\\
g_{t_1}(\lambda)g_{t_2}(\lambda)=g_{t_1+t_2}(\lambda)\text{ for }t_1t_2>0,\\
g_{t_1}(\lambda)g_{t_2}(\lambda)=0\text{ for }t_1t_2<0,\\
\frac{d}{dt}g_{t}(\lambda)=\lambda g_{t}(\lambda).\qed
\end{gather*}
\renewcommand\qed{}
\end{proof}

The main property of Green's function is described in the following theorem.
\begin{theorem}[{\rm\cite[Theorem 4.1, p.~81]{Daletskii-Krein:eng}}]\label{t:Green}
Equation~\eqref{e:x'=Ax+f} has a unique bounded on $\mathbb R$ continuously differentiable solution $x$ for any bounded continuous
function $f$ if and only if the spectrum $\sigma(A)$ does not intersect the
imaginary axis.
This solution possesses the representation
\begin{equation*}
x(t)=\int_{-\infty}^\infty \mathcal G(t-s)f(s)\,ds,
\end{equation*}
where $\mathcal G$ is Green's function~\eqref{e:Green's function} of equation~\eqref{e:x'=Ax+f}.
\end{theorem}

\begin{remark}\label{r:A-norm est}
We note that knowing an estimate of the function $t\mapsto\Vert\mathcal G(t)\Vert$ is an important information in the freezing method for equations with slowly varying coefficients~\cite[\S~10.2]{Bylov-Vinograd:rus1},~\cite[\S~7.4]{Henry81:eng},~\cite[ch.~10, \S~3]{Levitan-Zhikov82:eng},~\cite{Baskakov93:eng,Behncke-Hinton-Remling,Kuznetsova90:eng,Kuznetsova03:eng,Potzsche,Robinson,Xiao}. See also references therein.
\end{remark}

Below we assume that $A$ is a fixed complex $N\times N$-matrix and its spectrum does not intersect the imaginary axis. We denote by $\mu_1,\dots,\mu_k$ the roots of the characteristic polynomial that lie in the open right half-plane $\Real\mu>0$ counted according to their multiplicities; and we denote by $\nu_1,\dots,\nu_m$ the roots of the characteristic polynomial that lie in the open left half-plane $\Real\nu<0$ counted according to their multiplicities. Thus, $k+m=N$.

\begin{proposition}\label{p:mu-nu:2}
Let an analytic function $f$ be identically zero in the open right half-plane $\Real\mu>0$ {\rm(}an example of such a function is the function $\exp^+_t${\rm)}. Then
\begin{equation*}
f[\mu_1,\dots,\mu_k;\nu_{1},\dots,\nu_m]=\tilde f[\nu_{1},\dots,\nu_m],
\end{equation*}
where
\begin{equation*}
\tilde f(z)=\frac{f(z)}{\prod_{i=1}^k(z-\mu_i)}.
\end{equation*}
\end{proposition}
\begin{proof}
Suppose that all multiplicities equal 1.
By Proposition~\ref{p:repr of Delta}
\begin{align*}
f[\mu_1,\dots,\mu_k;\nu_{1},\dots,\nu_m]&=\sum_{q=1}^m\frac{f(\nu_q)}
{\prod\limits_{i=1}^k(\nu_q-\mu_i)\prod\limits_{\substack{j=1\\j\neq q}}^m(\nu_q-\nu_j)}\\
&=\sum_{q=1}^m\frac{\frac{f(\nu_q)}
{\prod\limits_{i=1}^k(\nu_q-\mu_i)}}{\prod\limits_{\substack{j=1\\j\neq q}}^m(\nu_q-\nu_j)}=\tilde f[\nu_{1},\dots,\nu_m].
\end{align*}
Since divided differences continuously depend on their arguments (Corollary~\ref{c:f[] is continuous}), the case of multiple points of interpolation is obtained by a passage to the limit.
\end{proof}

\begin{theorem}\label{t:p_t}
Let us arrange the roots of the characteristic polynomial in the following order{\rm:}
\begin{equation}\label{e:order}
\mu_1,\dots,\mu_k;\,\nu_1,\dots,\nu_m.
\end{equation}
Then the Newton interpolating polynomial $p_t^+$ of the function $\exp^+_t$ takes the form
\begin{equation}\label{e:p+}
p_t^+(z)=(z-\mu_1)\dots(z-\mu_k)q_t^+(z),
\end{equation}
where
\begin{equation*}
q_t^+(z)=\widetilde{\exp^+_t}[\nu_{1}]+\dots+\widetilde{\exp^+_t}[\nu_{1},\dots,\nu_m](z-\nu_{1})\dots(z-\nu_{m-1})
\end{equation*}
is the interpolating polynomial of the function
\begin{equation*}
\widetilde{\exp^+_t}(z)=\frac{\exp^+_t(z)}{\prod_{i=1}^k(z-\mu_i)}
\end{equation*}
with respect to the points $\nu_1,\dots,\nu_m$.
The interpolating polynomial $p_t^-$ of the function $\exp^-_t$ can be represented in the form
\begin{equation}\label{e:p-}
p_t^-(z)=(z-\nu_1)\dots(z-\nu_m)q_t^-(z),
\end{equation}
where
\begin{equation*}
q_t^-(z)=\widetilde{\exp^-_t}[\mu_{1}]+\dots+\widetilde{\exp^-_t}[\mu_{1},\dots,\mu_k](z-\mu_{1})\dots(z-\mu_{k-1})
\end{equation*}
is the interpolating polynomial of the function
\begin{equation*}
\widetilde{\exp^-_t}(z)=\frac{\exp^-_t(z)}{\prod_{j=1}^m(z-\nu_i)}
\end{equation*}
with respect to the points $\mu_1,\dots,\mu_k$.
\end{theorem}
\begin{proof}
We observe that $\exp^+_t(\mu_i)=0$, $i=1,\dots,k$. Therefore
\begin{equation*}
\exp^+_t[\mu_1]=\dots=\exp^+_t[\mu_1,\dots,\mu_k]=0.
\end{equation*}
Now from Proposition~\ref{p:Newton poly} it follows that
\begin{align*}
p_t^+(z)&=\exp^+_t[\mu_1,\dots,\mu_k;\nu_{1}](z-\mu_1)\dots(z-\mu_k)+\dots\\
&+\exp^+_t[\mu_1,\dots,\mu_k;\nu_{1},\dots,\nu_m](z-\mu_1)\dots(z-\mu_k)(z-\nu_{1})\dots(z-\nu_{m-1}).
\end{align*}
It remains to apply Proposition~\ref{p:mu-nu:2}.
\end{proof}

We set
\begin{align*}
\pi^+(\lambda)&=\begin{cases}
1, & \text{if $\Real\lambda<0$},\\
0, & \text{if $\Real\lambda>0$},\end{cases}\qquad
\pi^-(\lambda)=\begin{cases}
0, & \text{if $\Real\lambda<0$},\\
-1, & \text{if $\Real\lambda>0$}.\end{cases}
\end{align*}
Clearly, $\pi^+=\exp^+_{+0}=g_{+0}$ and $\pi^-=\exp^-_{-0}=g_{-0}$. Thus $P^+=\pi^+(A)$ and $P^-=\pi^-(A)$ are the spectral projectors that correspond to the subsets $\{\,\mu_1,\dots,\mu_k\,\}\subseteq\sigma(A)$ and $\{\,\nu_1,\dots,\nu_m\,\}\subseteq\sigma(A)$, respectively. Cf. Proposition~\ref{p:G-properties}.

In passing, we note the following corollary.

\begin{corollary}\label{c:pi_}
Let us arrange the roots of the characteristic polynomial in order~\eqref{e:order}.
Then the Newton interpolating polynomial $s^+$ of the function $\pi^+$ takes the form
\begin{equation*}
s^+(z)=(z-\mu_1)\dots(z-\mu_k)r^+(z),
\end{equation*}
where
\begin{equation*}
r^+(z)=\widetilde{\pi^+}[\nu_{1}]+\dots+\widetilde{\pi^+}[\nu_{1},\dots,\nu_m](z-\nu_{1})\dots(z-\nu_{m-1})
\end{equation*}
is the interpolating polynomial of the function
\begin{equation*}
\widetilde{\pi^+}(z)=\frac{1}{\prod_{i=1}^k(z-\mu_i)}
\end{equation*}
with respect to the points $\nu_1,\dots,\nu_m$.
The interpolating polynomial $s^-$ of the function $\pi^-$ can be represented in the form
\begin{equation*}
s^-(z)=(z-\nu_1)\dots(z-\nu_m)r^-(z),
\end{equation*}
where
\begin{equation*}
r^-(z)=\widetilde{\pi^-}[\mu_{1}]+\dots+\widetilde{\pi^-}[\mu_{1},\dots,\mu_k](z-\mu_{1})\dots(z-\mu_{k-1})
\end{equation*}
is the interpolating polynomial of the function
\begin{equation*}
\widetilde{\pi^-}(z)=\frac{1}{\prod_{j=1}^m(z-\nu_i)}
\end{equation*}
with respect to the points $\mu_1,\dots,\mu_k$.
\end{corollary}

\section{Sensitivity}\label{s:sensitivity}
Let $X$ be a Banach space. We denote by $\mathoo B(X)$ the algebra of all bounded linear operators $A:\,X\to X$. If $X=\mathbb C^N$, we identify $\mathoo B(X)$ with the algebra of $N\times N$-matrices.

Let $t\neq0$ be fixed.
The \emph{{\rm(}Fr\'echet{\rm)} differential} of the mapping $A\mapsto g_t(A)$
is defined~\cite[ch.~8]{Dieudonne:eng} as a linear mapping
$dg_t(\cdot,A):\;\mathoo B(\mathbb C^N)\to\mathoo B(\mathbb C^N)$ depending on the matrix parameter $A$ that possesses the property
\begin{equation}\label{e:dg_t}
g_t(A+\Delta A)-g_t(A)=dg_t(\Delta A,A)+o(\Vert\Delta A\Vert).
\end{equation}
The mapping $dg_t(\cdot,A)$ characterizes~\cite[Theorem 3.1]{Higham08} the influence on $g_t(A)$ of small perturbations in $A$. Such errors are unavoidable in the sense that they can not be diminished, regardless of which method is used to compute $\mathcal G(t)=g_t(A)$. The main characteristic of the magnitude of such errors is the number $\varkappa(g_t,A)=\Vert dg_t(\cdot,A)\Vert$. It is called the \emph{condition number} of the function $g_t$ at the point $A$.

\begin{proposition}[{\rm\cite[Theorem 3.9]{Higham08},~\cite[Theorems 67 and 71]{Kurbatov-Kurbatova-Oreshina}}]\label{p:dg_t=boxtimes}
Let the spectrum of the matrix $A$ do not intersect the imaginary axis. Then
\begin{align*}
dg_t(\Delta A,A)&=\frac1{2\pi i}\int_\Gamma g_t(\lambda)(\lambda\mathbf1-A)^{-1}\;\Delta A\;(\lambda\mathbf1-A)^{-1}\, d\lambda\\
&=\frac1{(2\pi
i)^2}\int_{\Gamma_1}\int_{\Gamma_2}g_t[\lambda,\mu](\lambda\mathbf1-A)^{-1}\;\Delta A\;(\mu\mathbf1-A)^{-1}\,d\mu\,d\lambda,
\end{align*}
where the contours $\Gamma$, $\Gamma_1$, and $\Gamma_2$ surround the spectrum of $A$.
The spectrum of the mapping $dg_t(\cdot,A)$ is the set
\begin{equation}\label{e:sigma(g_t[])}
\sigma\bigl[dg_t(\cdot,A)\bigr]=
\bigl\{\,g_t[\lambda,\mu]:\,\lambda,\mu\in\sigma(A)\,\bigr\}.
\end{equation}
\end{proposition}

We note that
\begin{equation*}
g_t[\lambda,\mu]=\begin{cases}
\exp^+_t[\lambda,\mu], & \text{if $t>0$},\\
-\exp^-_t[\lambda,\mu], & \text{if $t<0$},
\end{cases}
\end{equation*}
where
\begin{align*}
\exp^+_t[\lambda,\mu]&=\begin{cases}
\frac{e^{\lambda t}-e^{\mu t}}{\lambda-\mu}, & \text{if $\Real\lambda<0$ and $\Real\mu<0$},\\
te^{\lambda t}, & \text{if $\Real\lambda<0$ and $\lambda=\mu$},\\
\frac{e^{\lambda t}}{\lambda-\mu}, & \text{if $\Real\lambda<0$ and $\Real\mu>0$},\\
-\frac{e^{\mu t}}{\lambda-\mu}, & \text{if $\Real\lambda>0$ and $\Real\mu<0$},\\
0, & \text{if $\Real\lambda>0$ and $\Real\mu>0$},
\end{cases}
\end{align*}
\begin{align*}
\exp^-_t[\lambda,\mu]&=\begin{cases}
0, & \text{if $\Real\lambda<0$ and $\Real\mu<0$},\\
-\frac{e^{\mu t}}{\lambda-\mu}, & \text{if $\Real\lambda<0$ and $\Real\mu>0$},\\
\frac{e^{\lambda t}}{\lambda-\mu}, & \text{if $\Real\lambda>0$ and $\Real\mu<0$},\\
\frac{e^{\lambda t}-e^{\mu t}}{\lambda-\mu}, & \text{if $\Real\lambda>0$ and $\Real\mu>0$},\\
te^{\lambda t}, & \text{if $\Real\lambda>0$ and $\lambda=\mu$}.
\end{cases}
\end{align*}

The following theorem is an analogue of the representation for the differential of the matrix exponential found in~\cite[formula (1.8)]{Karplus-Schwinger}, see also~\cite[ch.~10, \S~14]{Bellman69},~\cite[formula (10.15)]{Higham08}, \cite[example 2]{Kenney-Laub89},~\cite{van_Loan77}.

\begin{theorem}\label{t:Karplus-Schwinger}
The following representation holds{\rm:}
\begin{equation*}
dg_t(\Delta A,A)=\int_{-\infty}^\infty g_s(A)\;\Delta A\;g_{t-s}(A)\,ds.
\end{equation*}
\end{theorem}
\begin{proof}
By Proposition~\ref{p:dg_t=boxtimes} (see also~\cite[Theorem 32]{Kurbatov-Kurbatova-Oreshina}) it is enough to establish the identity
\begin{equation*}
g_t[\lambda,\mu]=\int_{-\infty}^\infty g_s(\lambda)g_{t-s}(\mu)\,ds.
\end{equation*}

By the definition, we have
\begin{align*}
g_t(\lambda)&=\begin{cases}
0, & \text{if $t<0$},\\
e^{\lambda t}, & \text{if $t>0$},
\end{cases}&\text{for }\lambda&<0,\\
g_t(\lambda)&=\begin{cases}
-e^{\lambda t}, & \text{if $t<0$},\\
0, & \text{if $t>0$}
\end{cases}&\text{for }\lambda&>0.
\end{align*}

We consider several cases.
Suppose that $t>0$.

Suppose that $\Real\lambda<0$ and $\Real\mu<0$, $\lambda\neq\mu$. Then we have
\begin{equation*}
\int_{-\infty}^\infty g_s(\lambda)g_{t-s}(\mu)\,ds=\int_{(0,\infty)\cap(-\infty,t)} e^{\lambda s}e^{\mu(t-s)}\,ds=
\frac{e^{\lambda t}-e^{\mu t}}{\lambda-\mu}=\exp^+_t[\lambda,\mu].
\end{equation*}
Similarly, for $\Real\lambda<0$ and $\Real\mu<0$, $\lambda=\mu$, we have
\begin{equation*}
\int_{-\infty}^\infty g_s(\lambda)g_{t-s}(\mu)\,ds
=\int_{(0,\infty)\cap(-\infty,t)} e^{\lambda s}e^{\lambda(t-s)}\,ds=te^{\lambda t}=\exp^+_t[\lambda,\lambda].
\end{equation*}

Suppose that $\Real\lambda<0$ and $\Real\mu>0$. Then
we have
\begin{align*}
\int_{-\infty}^\infty g_s(\lambda)g_{t-s}(\mu)\,ds&=\int_{(0,\infty)\cap(t,\infty)}e^{\lambda s}\bigl(-e^{\mu(t-s)}\bigr)\,ds=
\frac{e^{\lambda t}}{\lambda-\mu}=\exp^+_t[\lambda,\mu].
\end{align*}

Suppose that $\Real\lambda>0$ and $\Real\mu<0$. Then we have
\begin{align*}
\int_{-\infty}^\infty g_s(\lambda)g_{t-s}(\mu)\,ds&=\int_{(-\infty,0)\cap(-\infty,t)}\bigl(-e^{\lambda s}\bigr)e^{\mu(t-s)}\,ds=
-\frac{e^{\mu t}}{\lambda-\mu}=\exp^+_t[\lambda,\mu].
\end{align*}

Suppose that $\Real\lambda>0$ and $\Real\mu>0$. Then
we have
\begin{align*}
\int_{-\infty}^\infty g_s(\lambda)g_{t-s}(\mu)\,ds&=\int_{(-\infty,0)\cap(t,+\infty)}\bigl(-e^{\lambda s}\bigr)\bigl(-e^{\mu(t-s)}\bigr)\,ds=0=\exp^+_t[\lambda,\mu].
\end{align*}

The case $t<0$ is considered analogously.
\end{proof}

\begin{corollary}\label{c:est of dg_t}
The following estimates hold{\rm:}
\begin{align}
\Vert dg_t(\cdot,A)\Vert&\le\int_{-\infty}^\infty\Vert g_s(A)\Vert\cdot\Vert g_{t-s}(A)\Vert \,ds,\label{e:est}\\
\int_{-\infty}^\infty\Vert dg_t(\cdot,A)\Vert\,dt&\le\Bigl(\int_{-\infty}^\infty\Vert g_s(A)\Vert\,ds\Bigr)^2.\notag
\end{align}
\end{corollary}
\begin{proof}
From Theorem~\ref{t:Karplus-Schwinger} we have
\begin{equation*}
\Vert dg_t(\Delta A,A)\Vert\le\int_{-\infty}^\infty\Vert g_s(A)\Vert\cdot\Vert\Delta A\Vert\cdot\Vert g_{t-s}(A)\Vert \,ds,
\end{equation*}
which implies~\eqref{e:est}. The second inequality follows from the first one and the properties of convolution.
\end{proof}

\section{The algorithm}\label{s:algorithm}
In this Section we describe an algorithm based on Theorem~\ref{t:p_t} for computing Green's function of bounded solutions problem. The algorithm is especially convenient for symbolic calculations.

0. Given a complex square matrix $A$. Let its size be $N\times N$. If $N$ is small ($N\le10$), Green's function can be calculated symbolically, i.~e. presented as an expression that depends on $t$. If $N$ is large, we choose a point $t\neq0$ and perform the following calculations for the specific value of $t$. If $N$ is large, but Green's function is needed in the form of an expression depending on $t$, we calculate $\mathcal G(t)$ at several points $t$ and then use an interpolation.

1. We calculate eigenvalues $\lambda_1,\dots,\lambda_N$ of the matrix $A$ counted according to their multiplicities. For example, it can be done by means of the Schur algorithm (see, e.g.,~\cite[ch.~7]{Golub-Van_Loan96}). (We note that the employment of the Schur triangular form may accelerate the subsequent substitution of the matrix $A$ into polynomials.)

2. We divide eigenvalues into two groups: $\mu_1,\dots,\mu_k$ lie in the open right half-plane $\Real\mu>0$ and $\nu_1,\dots,\nu_m$ lie in the open left half-plane $\Real\nu<0$. At the same time, we verify if any eigenvalue lies on the imaginary axis.

3. We calculate divided differences of the functions
\begin{align*}
\widetilde{\exp^+_t}(z)&=\frac{\exp^+_t(z)}{\prod_{i=1}^k(z-\mu_i)}
=\frac{e^{z t}}{\prod_{i=1}^k(z-\mu_i)}
,&t&>0,\;\Real z<0,\\
\widetilde{\exp^-_t}(z)&=\frac{\exp^-_t(z)}{\prod_{j=1}^m(z-\nu_i)}
=\frac{e^{z t}}{\prod_{j=1}^m(z-\nu_i)},&t&<0,\;\Real z>0.
\end{align*}
If one renumbers $\mu_1,\dots,\mu_k$ so that $\Real\mu_1\ge\Real\mu_2\ge\dots\ge\Real\mu_k$ and renumbers $\nu_1,\dots,\nu_m$ so that $\Real\nu_1\le\Real\nu_2\le\dots\le\Real\nu_m$,
then the first divided differences are significantly reduced.

4. According to Proposition~\ref{p:f->p} and Theorem~\ref{t:p_t} we calculate Green's function $\mathcal G$ by the formulas
\begin{align*}
\mathcal G(t)&=p_t^+(A),&t&>0, \\
\mathcal G(t)&=-p_t^-(A),&t&<0,
\end{align*}
where $p_t^+$ and $p_t^-$ are defined by formulas~\eqref{e:p+} and~\eqref{e:p-}.
It is convenient to calculate $p_t^\pm(A)$ by the rules
\begin{align*}
p_t^+(A)&=R_1\prod_{i=1}^k(A-\mu_i\mathbf1),\\
p_t^-(A)&=S_1\prod_{j=1}^m(A-\nu_j\mathbf1),
\end{align*}
where $R_1=q_t^+(A)$ and $S_1=q_t^-(A)$ are calculated according to the Horner algorithm:
\begin{align*}
R_m&=\widetilde{\exp^+_t}[\nu_{1},\dots,\nu_m]\mathbf1,\\
R_{j-1}&=(A-\nu_{j-1}\mathbf1)R_j+\widetilde{\exp^+_t}[\nu_{1},\dots,\nu_{j-1}]\mathbf1, &j&=m,m-1,\dots,2;\\
S_k&=\widetilde{\exp^-_t}[\nu_{1},\dots,\nu_m]\mathbf1,\\
S_{i-1}&=(A-\mu_{i-1}\mathbf1)S_i+\widetilde{\exp^-_t}[\mu_{1},\dots,\mu_{i-1}]\mathbf1, &i&=k,k-1,\dots,2.
\end{align*}
We stress that the products $\prod_{i=1}^k(A-\mu_i\mathbf1)$ and $\prod_{j=1}^m(A-\nu_j\mathbf1)$ do not contain $t$. Therefore it is enough to calculate them only once.

5. As a verification of the plausibility of the results obtained, we propose to verify the identities from Proposition~\ref{p:G-properties}.

\begin{remark}\label{r:kk}
If there are close eigenvalues, large rounding errors can occur, see the discussion of this phenomenon in~\cite{Moler-Van_Loan03}. In such a case, the method of calculating divided differences proposed in~\cite{Kurbatov-Kurbatova-LMA-2016} can be applied.
\end{remark}

\begin{remark}\label{r:P pm}
One can propose a similar algorithm based on Corollary~\ref{c:pi_} for the calculation of the spectral projectors $P^\pm$.
\end{remark}

\section{Numerical experiments}\label{s:experiments}
For numerical experiments we took matrices $A$ consisting of random (uniformly distributed) complex numbers from the rectangle $[-1,1]\times[-i,i]\subset\mathbb C$.
Numerical experiments showed that the Algorithm is reliable for $N\le40$. If $N>50$, the relative errors are about 20\%.
The condition number is usually less than $10^3$.

As noted above, if $N\le 10$, it is possible to obtain Green's function in symbolic from, i.~e. as an expression depending on $t$. In this case, the identities from Proposition~\ref{p:G-properties} were fulfilled with the accuracy $10^{-10}$.

\begin{figure}[htb]
\begin{center}
\includegraphics[width=.94\textwidth]{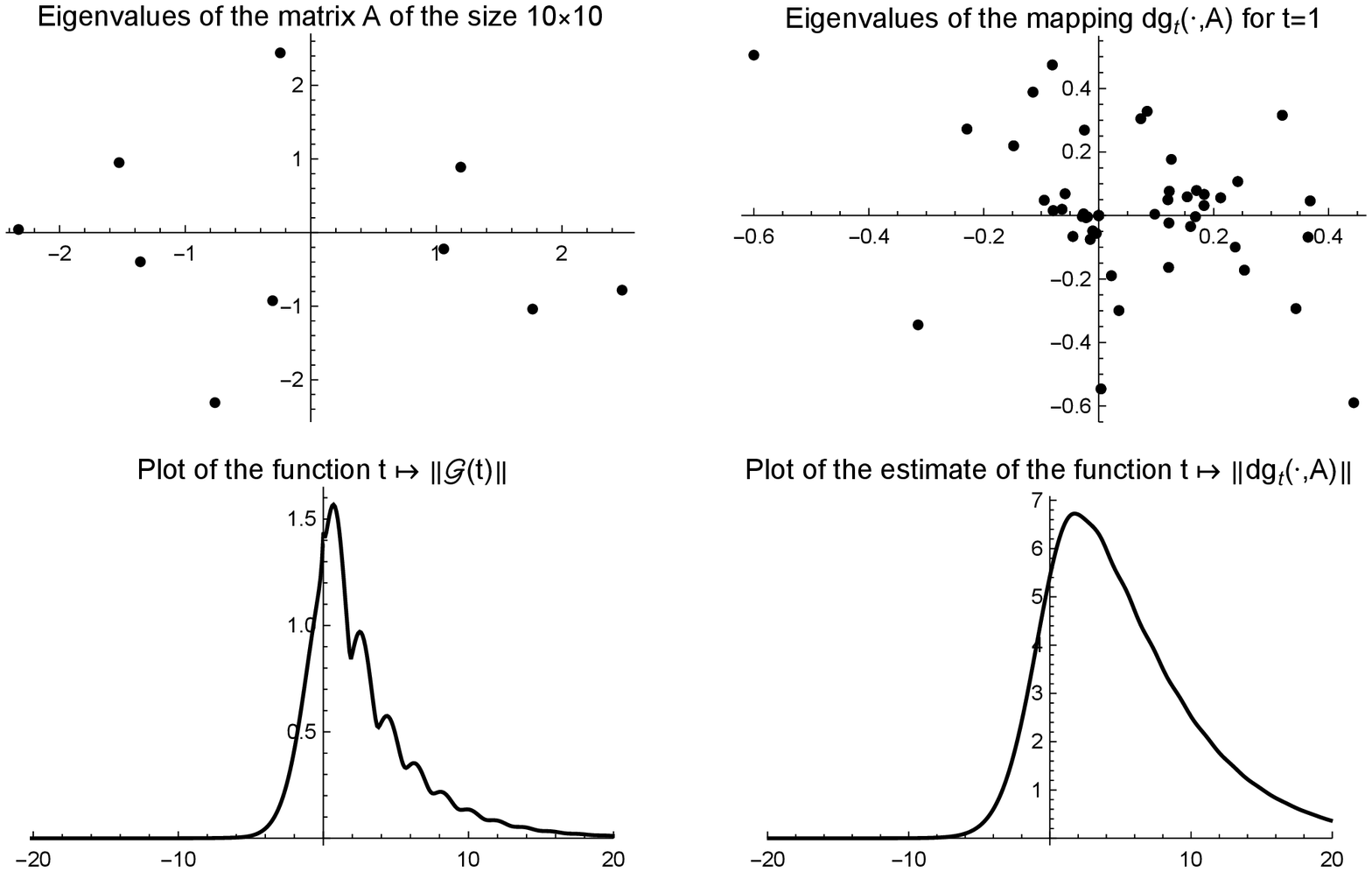}
\caption{Numerical experiment with a $10\times10$-matrix} \label{f:green10}
\end{center}
\bigskip
\begin{center}
\includegraphics[width=.94\textwidth]{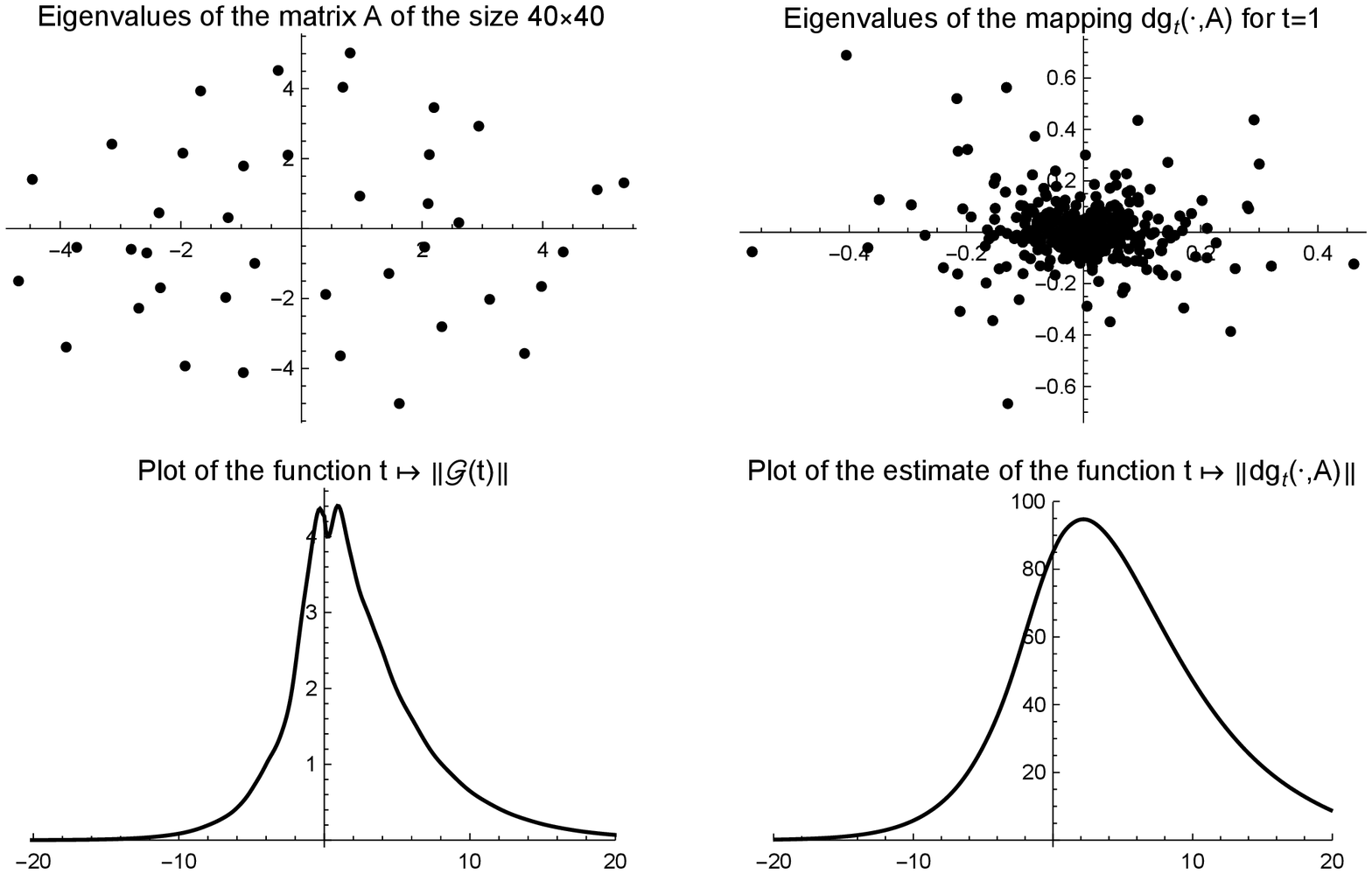}
\caption{Numerical experiment with a $40\times40$-matrix} \label{f:green40}
\end{center}
\end{figure}

We present results of two numerical experiments on fig.~\ref{f:green10}--\ref{f:green40}.
The eigenvalues of $dg_t(\cdot,A)$ were calculated according to formula~\eqref{e:sigma(g_t[])}. We took the Euclidian norm in $\mathbb C^N$ and associated norms in $\mathoo B(\mathbb C^N)$ and $\mathoo B\bigl(\mathoo B(\mathbb C^N)\bigr)$; we calculated the norm of a matrix as its largest singular number. Estimates of the norms of $dg_t(\cdot,A)$ were calculated by formula~\eqref{e:est}.

\end{document}